\documentclass[12pt]{amsart}
\usepackage[english]{babel}
\usepackage[all]{xy}
 \usepackage[latin1]{inputenc}
  \usepackage[T1]{fontenc}
  \usepackage{amsmath,amssymb}
  \usepackage{amsmath}
  \usepackage{amsthm}
\usepackage{amssymb}
\usepackage{mathrsfs}
\usepackage{enumerate}
\usepackage{tikz}
\usetikzlibrary{shapes.geometric}
\usepackage{physics}
\usepackage[notcite, final, notref]{showkeys}
\usepackage{dsfont}
\newtheorem{theorem}{Theorem}[section]

\newtheorem{proposition}[theorem]{Proposition}

\newtheorem{definition}[theorem]{Definition}

\theoremstyle{definition}
\newtheorem{remark}[theorem]{Remark}

\newtheorem{example}[theorem]{Example}

\newcommand{\w}{\omega}

\usepackage{tikz}
\usetikzlibrary{arrows,shapes,positioning,shadows,trees,matrix}
\usepackage{xcolor}

\title[]{Similarity metrics, metrics, and conditionally negative definite functions}

\author[D. Alpay]{Daniel Alpay}
\address{(DA) Schmid College of Science and Technology \\
Chapman University\\
One University Drive
Orange, California 92866\\
USA}
\email{alpay@chapman.edu}

\author[L. Mayats-Alpay]{Liora Mayats-Alpay}
\address{(LMA) Fowler School of Engineering and Schmid College of Science and Technology\\
Chapman University\\
One University Drive
Orange, California 92866\\
USA}
\email{mayatsalpay@chapman.edu}

\thanks{D. Alpay thanks the Foster G. and Mary McGaw Professorship in
  Mathematical Sciences, which supported his research.}

\begin{document}
\maketitle

\begin{abstract}
  We give an example of a similarity metric which is not positive definite, and present a general theorem which provides a large family of similarity metrics which are positive
  definite.
\end{abstract}

\tableofcontents
\section{Introduction}
\setcounter{equation}{0}
{\bf Prologue:}
Among the numerous deep mathematical tools used in machine learning, an important role is played by metric and metric spaces (see Definition \ref{metric-space})
and positive definite functions (see Definition \ref{pdfunction}), for instance in the $k$-neighbors algorithm (KNN) or in the theory of support vector machines (SVM).
These two notions are connected by different threads. For example, one can associate to a positive definite function on a set $X$ a metric on $X$
(possibly after moding out points corresponding to the same kernel function). See for instance \cite{MR4302453} and Theorem \ref{dkdk} below.
Let now $(X,D)$ be a metric space. When the metric is normalized, that is, takes values in $[0,1]$, the function
\begin{equation}
  \label{normalized}
  s(x,y)=1-D(x,y)
\end{equation}
is used to study the similarity of $x$ and $y$, and bears various names, such as similarity metric (note that a similarity metric is never
a metric). See \cite{MR2522441}.\smallskip

The function $s(x,y)$ satisfies the following properties, which follow directly from the fact that $D(x,y)$ is a normalized metric.
\begin{eqnarray}
  \label{paris-1}
  s(x,y)&
          \in&[0,1],\\
  s(x,y)&=&s(y,x),\\
  s(x,x)&=&1,\\
  s(x,z)+s(z,y)&\le & s(x,y) +1\\
  s(x,y)=1&\iff &x=y.
                  \label{paris-6}
  \end{eqnarray}

Conversely, any function on $X\times X$ satisfying the above properties will define a normalized metric via \eqref{normalized}. See \cite[Corollary 1 p. 2367]{MR2522441}.\\
  
{\bf Similarity functions and Schoenberg's theorem:}
When $d(x,y)$ is a metric on the space $X$, the function
\begin{equation}
  \label{pre-shafer}
D(x,y)= 1-e^{-d(x,y)}
\end{equation}
is a normalized metric on $X$ (see \cite[Lemma 12 p. 2370]{MR2522441}). Equation \eqref{normalized} becomes
\begin{equation}
  s(x,y)=e^{-d(x,y)}.   
  \end{equation}
A theorem of Schoenberg   gives a characterization of the functions $d(x,y)$ 
for which $e^{-td(x,y)}$ is positive definite on $X$ for every $t>0$. These are exactly the conditionally negative functions (see Theorem \ref{Schoenber} below).
This suggest possible connections between positive definite functions and similarity metrics. Note that $d(x,y)$ in Schoenberg's theorem need not be a metric, but
in the real-valued setting, a conditionally negative functions is the square of a (pseudo-)distance. See \cite[p. 85]{MR86b:43001}. We also note that In Schoenberg's theorem, the function $e^{-td(x,y)}$ can be replaced by $\frac{1}{1+td(x,y)}$ (and by a host of other functions
obtained using Laplace transform; see \cite[Theorem 2.3 p.  75]{MR86b:43001}).\\
\\


{\bf Similarity functions and evidence theory:} Equality \eqref{pre-shafer}
appears also in evidence theory. There $d$ is a weight of evidence, and $D$ is the degree of support of the evidence. See \cite[pp. 8-9]{shafer1976mathematical}.\\

{\bf Some questions:} Motivated by Schoenberg's theorem and in part by the work  \cite{song2015distance}, the following questions seem natural to raise
(there is some overlap between the various questions): \\
$(a)$ Given a normalized metric $D(x,y)$ on a set $X$ when is $1-D(x,y)$ positive definite on $X$?\\
Conversely,\\
$(b)$ Given a positive definite function $s(x,y)$ taking values in $[0,1]$, when is it a similarity metric, that is when is $1-s(x,y)$ a (normalized) metric?\\

{\bf Outline:} The paper consists of five sections besides the introduction. In Section 2 we review a few facts on metrics and Schoenberg's theorem. In Section
3 we discuss similarity metrics. The new results are presented in Section 4 and 5. In Section 4 we give an example of a similarity metric which is not
positive definite. In Section 5 we present a large class of positive definite similarity metrics. The last section is devoted to some more examples of positive
definite similarity metrics.

\section{Metric, positive definite and negative definite functions}
\setcounter{equation}{0}

We begin by recalling the definition of a metric space since the notion is also used in the present section.

\begin{definition}
  \label{metric-space}
Let $X$ be a set. The map $d(x,y)$ from $X\times X$ into $[0,\infty)$ is called a metric (or a distance) if the following three conditions hold for all $x,y,z\in X$:
  \begin{eqnarray}    \label{sym5}
    d(x,y)&=& 0\quad\iff\quad x=y\\
    \label{sym6}
    d(x,y)&=&d(y,x)\\
    d(x,y)&\le&d(x,z)+d(y,z).
                \label{triangu}
\end{eqnarray}
  \end{definition}

  We now turn to the definitions of positive definite and negative definite kernels. The latter are also called {\sl conditionally negative definite kernels}
  and are {\sl not} the opposite of positive definite kernels. We will also use the term {\sl function} rather than  {\sl kernel}, although this terminology
  refers in general to a smaller class of kernels.
  
\begin{definition}
  \label{pdfunction}
  Let $X$ be a set. The map $k(x,y)$ from $X\times X$ into $\mathbb C$ is called positive definite if the following condition holds: For every $N\in\mathbb N$, every  choice of $x_1,\ldots, x_N\in X$ and every choice of
  $c_1,\ldots, c_N\in\mathbb C$ it holds that
  \[
\sum_{j,k=1}^N\overline{c_k}k(x_k,x_j)c_j\ge 0.
    \]
\end{definition}

\begin{definition}
Let $X$ be a set. The map $k(x,y)$ from $X\times X$ into $\mathbb C$ is called negative definite if the following condition holds: For every $N\in\mathbb N$, every  choice of $x_1,\ldots, x_N\in E$ and every choice of
$c_1,\ldots, c_N\in\mathbb C$ such that
\[
\sum_{j=1}^Nc_j=0,
  \]
it holds that
\begin{equation}
  \label{negsum}
\sum_{j,k=1}^N\overline{c_k}k(x_k,x_j)c_j\le 0.
    \end{equation}
\end{definition}

Connections between the two notions are presented in \cite{MR86b:43001}, where one can find in particular the following result, originally due to Schoenberg (see \cite[Theorem 4 and Theorem 5]{MR1503439}):

\begin{theorem} (\cite[Theorem 2.2 p. 74]{MR86b:43001})
  \label{Schoenber}
  Let  $X$ be a set and let $d(x,y)$ be a complex-valued function defined on $X\times X$. Then, $d(x,y)$ is negative definite if and only if $e^{-td(x,y)}$ is positive definite
  for all $t>0$.
  \end{theorem}

  As a consequence of this theorem, and as a special case of  \cite[Theorem 2.3 p.  75]{MR86b:43001} we have

  \begin{theorem}
  \label{Schoenber1}
  Let  $X$ be a set and let $d(x,y)$ be a complex-valued function defined on $X\times X$. Then, $d(x,y)$ is negative definite if and only if $\frac{1}{1+td(x,y)}$
  is positive definite for all $t>0$.
  \end{theorem}

Indeed, we have
\[
  \frac{1}{1+td(x,y)}=\int_0^\infty e^{-utd(x,y)}e^{-u}du
\]
corresponding to $d\mu(u)=e^{-u}du$ in \cite[Theorem 2.3 p.  75]{MR86b:43001}.\smallskip

The following important connection between metrics and positive definite functions is used in the sequel.
\begin{theorem}
  Let $k(x,y)$ be a positive definite function on the set $X$, and assume that
  \begin{equation}
    \label{not-eq}
x\not=y\quad\longrightarrow\quad k(\cdot,x)\not\equiv k(\cdot, y).
\end{equation}
Then the formula
\begin{equation}
  \label{superform}
d_k(x,y)=\sqrt{k(x,x)+k(y,y)-2{\rm Re}\, k(x,y)}
\end{equation}
defines a metric on $X$.
\label{dkdk}
\end{theorem}
As an illustration we have:

\begin{proposition}
  \label{hurst1}
  Let $\mathcal H$ be a Hilbert space and let $a\in(0,1)$. The function
  \begin{equation}
    \label{hurst}
\|x\|^{2a}+\|y\|^{2a}-\|x-y\|^{2a}
\end{equation}
is positive definite on $\mathcal H$. The associated metric $d_k$ is equal to
\begin{equation}
  \label{xa}
d_k(x,y)=\|x-y\|^a
  \end{equation}
\end{proposition}

\begin{proof}
  Applying to $z=\|x-y\|^2$ the formula (see \cite[Corollary 2.10 p. 78]{MR86b:43001} and \cite[formula $(8)$ p. 526]{MR1501980},
  the latter being equivalent to the former via a change of variable)
  \[
    z^a=\frac{a}{\Gamma(1-a)}\int_0^\infty(1-e^{-tz})\frac{dt}{t^{a+1}},\quad a\in(0,1),\
  \]
  valid for any complex number $z$ with positive real part, we obtain
  \[
\|x-y\|^{2a}=\frac{a}{\Gamma(1-a)}\int_0^\infty(1-e^{-t\|x-y\|^2})\frac{dt}{t^{a+1}}.
\]
It follows that the function $\|x-y\|^{2a}$ is negative definite and hence, by \cite[Lemma 2.1 p. 74]{MR86b:43001}, that \eqref{hurst}
is positive definite on $\mathcal H$.\smallskip

Condition \eqref{not-eq} is readily verified and so is \eqref{xa}.
\end{proof}

\begin{remark} That \eqref{xa} is a metric could be obtained more directly, using the fact that, for $a\in(0,1)$,  $d^a$ is a metric when $d$ is a metric.
  \end{remark}

  We also remark that, for $\mathcal H=\mathbb C$, the kernel \eqref{hurst} is the covariance function of the fractional Brownian motion, with Hurst parameter $a$ (usually denoted by $H$).\\
  \section{Similarity metrics}
\setcounter{equation}{0}
In this section we discuss similarity metrics and their connections to metrics and positive definite functions. For a  different approach to similarity metric and related
counterexamples, we refer to \cite{MR3906938}.

\begin{definition} (see \cite[Definition 2 p. 2366]{MR2522441})
  Let $X$ be a set. A similarity metric is a real-valued map defined on $X\times X$ and satisfying the following five conditions for all $x,y,z\in X$:
  \begin{eqnarray}
    s(x,y)&=&s(y,x),\\
    s(x,x)&\ge&0,\\
    s(x,x)&\ge &s(x,y)\\
    s(x,y)+s(y,z)&\le &s(x,z)+s(y,y)\\
    s(x,x)=s(y,y)=s(x,y)&\iff& x=y.
                               \label{defchen}
    \end{eqnarray}
    \end{definition}

as opposed to the definition from \cite{song2015distance} (where $A,B$ and $C$ denotes intutionistic fuzzy sets). There,
the distance is on sets (more precisely, intuitionistic fuzzy sets) and the following conditions are added:
(we use letters $A,B, \ldots$ for the sets): $d(A,B)$ is normalized,
\begin{equation}
  0\le d(x,y)\le 1
\end{equation}
and the following condition, pertaining to inclusions of sets, holds:
\begin{equation}
  \label{songplus}
{\rm If}\quad A\subset B\subset C,\quad {\rm then}\quad d(A,B)\le d(A,C)\quad{\rm  and}\quad d(B,C)\le d(A,C).
\end{equation}
\begin{definition}
  \label{songdef}
  \begin{enumerate}
\item $0\le     S(A,B)\le 1$.
\item $S(A,B)=1$    if and only if $A=B$.
\item $S(A,B)=S(B,A)$.
  \item If $A\subset B\subset C$ then $S(A,B)\ge S(A,C)$ and $S(B,C)\ge S(A,C)$.
    \end{enumerate}
\end{definition}
These two definitions lead to two approaches to the notion of similarity.

\begin{theorem} (see \cite[Corollary 1 p.2367]{MR2522441})
  If $s$ is a positive normalized similarity metric in the sense of \eqref{paris-1}-\eqref{paris-6}, then $1-s(x,y)$ is a normalized metric, and conversely.
\end{theorem}



\section{A counterexample}
\setcounter{equation}{0}

We set $X$ to be the set of complex-valued piece-wise continuous functions defined on the real line, which are moreover bounded in absolute value, and on $X$ we define
\[
d(x,y)=\sup_{s\in\mathbb R}|x(s)-y(s)|
  \]
  Clearly $d$ is a metric on $X$, defined by the norm
  \begin{equation}
    \label{infty}
    \|x\|=\sup_{s\in\mathbb R}|x(s)|.
    \end{equation}
\begin{theorem}
  With $X$ and $d$ as above, the function $D(x,y)=1-e^{-d(x,y)}$ is a metric on $X$ and $1-D(x,y)$ is not positive definite on $X$
  \end{theorem}

  \begin{proof}
    We divide the proofs in a number of steps.\\

    STEP 1: {\sl $D(x,y)$ is a normalized metric.}\smallskip

    See \cite[Lemma 12 p. 2370]{MR2522441}.\\

    Consider now $d(x,y)$ restricted to the set $X_0$ of functions (see Figure 1)
   
    \begin{equation}
      \begin{split}
        x_1(s)&=\begin{cases}\, 1,\,\quad s\in[0,1],\\ \, 0,\,\quad{\rm otherwise.}\end{cases}\\
        x_2(s)&=\begin{cases}\, 1,\,\quad s\in[2,3],\\
          \, -1,\,\hspace{1mm} s\in[6,7],\\
          \, 0,\,\quad{\rm otherwise.}\end{cases}\\
                x_3(s)&=\begin{cases}\, 1,\,\quad s\in[4,5],\\
          \, -1,\,\hspace{1mm} s\in[2,3],\\
          \, 0,\,\quad{\rm otherwise.}\end{cases}\\
                x_4(s)&=\begin{cases}\, 1,\,\quad s\in[6,7],\\
          \, -1,\,\hspace{1mm} s\in[4,5],\\
          \, 0,\,\quad{\rm otherwise.}\end{cases}\\
        x_5(s)&=\, -1,\,\hspace{5mm} s\in[2,3].
      \end{split}
    \end{equation}

    \begin{figure}[h]
        \includegraphics[width=0.6\linewidth]{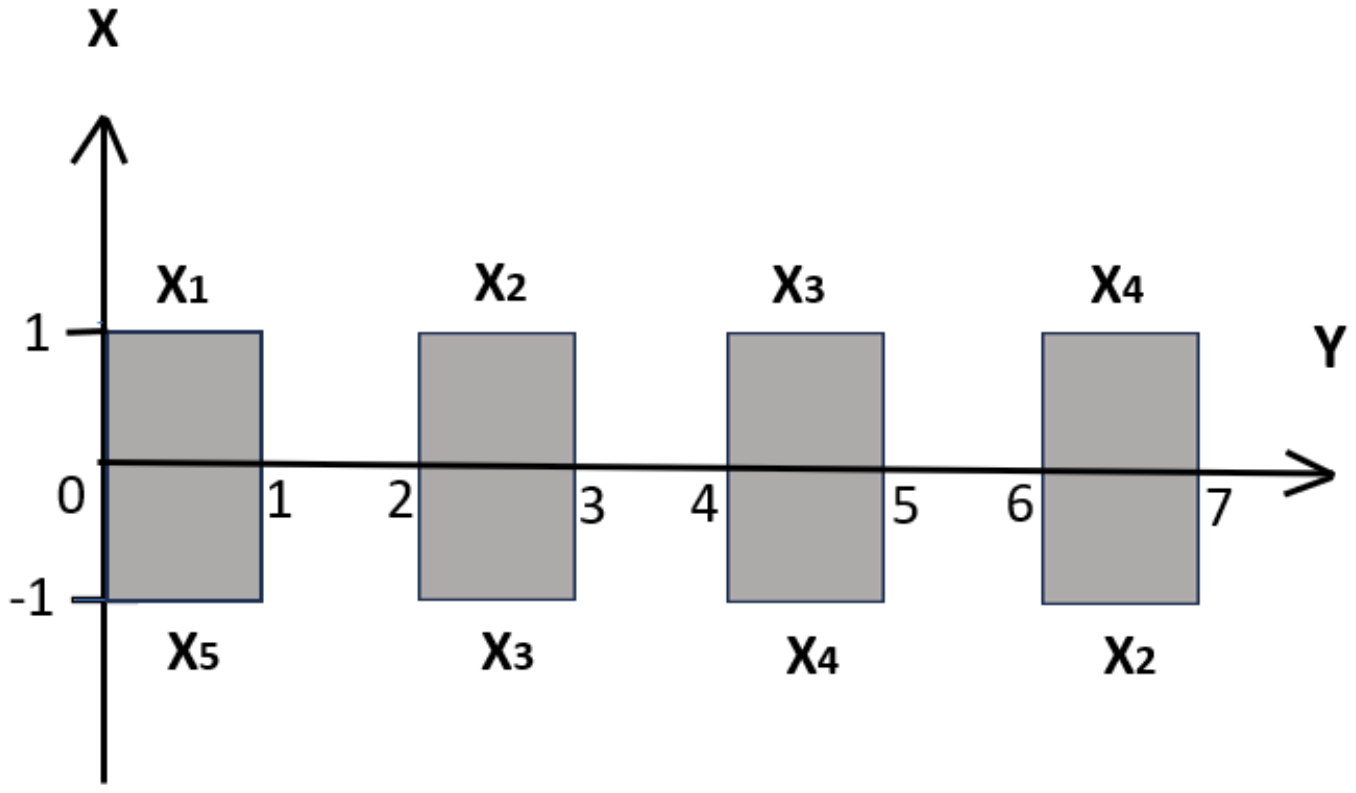}
        \vspace{-2cm}
        \caption{}
        \label{f1}
\end{figure}

    One has
        \[
      \begin{split}
        d(x_1,x_2)&=1\\
        d(x_1,x_3)&=1\\
        d(x_1,x_4)&=1\\
        d(x_1,x_5)&=2\\
                d(x_2,x_3)&=2\\
        d(x_2,x_4)&=2\\
        d(x_2,x_5)&=1\\
        d(x_3,x_4)&=2\\
        d(x_3,x_5)&=1\\
                        d(x_4,x_5)&=1.
        \end{split}
      \]
    The corresponding matrix $(d(x_j,x_k))_{j,k=1}^5$ is equal to
    \[
    \Delta=  \begin{pmatrix}0&1&1&1&2\\
        1&0&2&2&1\\
        1&2&0&2&1\\
        1&2&2&0&1\\
        2&1&1&1&0\end{pmatrix}.
    \]

    STEP 2:  {\sl $d(x,y)$ restricted to $X_0$ is not negative definite.}\smallskip

By \cite[Theorem 6.2.16 p. 80]{MR1460488} a necessary condition for $\Delta$ to be of  negative type is that it admits exactly one positive eigenvalue.
        A computation using the online matrix calculator https://matrixcalc.org/ shows that this matrix has $\lambda=-2$ as negative eigenvalue
    of multiplicity $3$ and two strictly positive eigenvalues, namely $\lambda=\pm \sqrt{7}+3$ with corresponding eigenvectors
    \[
      \begin{pmatrix}1\\ \frac{\sqrt{7}+1}{3}\\ \frac{\sqrt{7}+1}{3}\\ \frac{\sqrt{7}+1}{3}\\ 1      \end{pmatrix}\quad{\rm and}\quad      \begin{pmatrix}1\\ \frac{-\sqrt{7}+1}{3}\\ \frac{-\sqrt{7}+1}{3}\\ \frac{-\sqrt{7}+1}{3}\\ 1      \end{pmatrix}\quad{\rm respectively}.
    \]
Using the above cited theorem we conclude that $\Delta$ is not of negative type.\\
    
      STEP 3: {\sl We conclude using Schoenberg's theorem.}\smallskip

      By this theorem, there exists a $t>0$ for which $e^{-td(x,y)}$ is not positive definite. But in the case at hand,  $td(x,y)=d(tx,ty)$ for $t>0$ and so the condition is
      independent of $t>0$.
      \end{proof}

      We note that the norm \eqref{infty} is not strictly convex. We have
      \[
\|x_1-x_2\|+\|x_2-x_5\|=\|x_1-x_5\|
\]
but $x_1-x_2$ is not a multiple of $x_2-x_5$.\smallskip

We also note that the previous example can be written in terms of the following graphs, and the graph metric, called here $d_G$. Recall that the graph metric between
two nodes of a connected graph is the number of edges of the shortest path linking the two nodes. For the graph on the left in Figure 2 one has
    \begin{figure}[h]
        \includegraphics[width=0.6\linewidth]{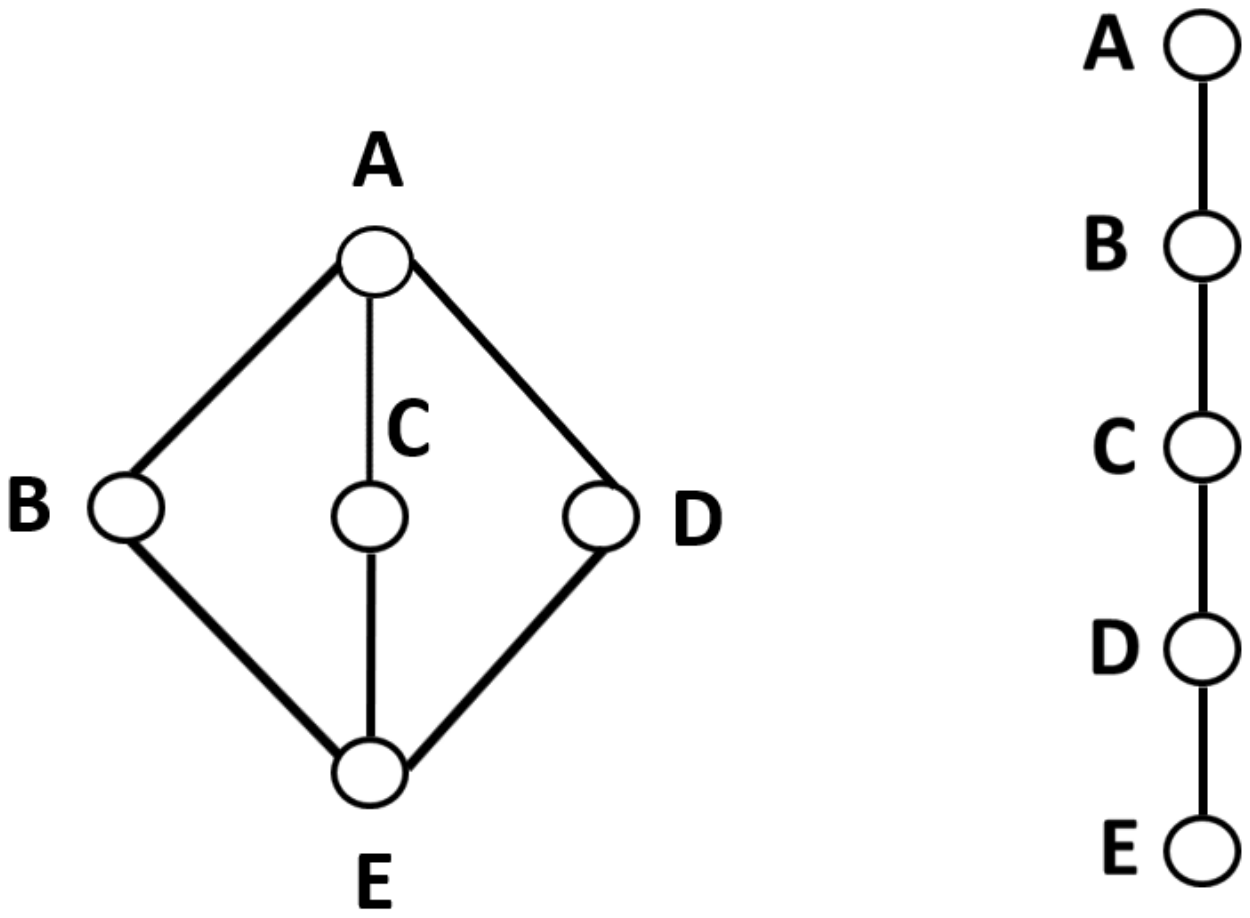}
        \vspace{-2cm}
\caption{}
\end{figure}

   \[
      \begin{split}
        d_G(A,B)&=1\\
        d_G(A,C)&=1\\
        d_G(A,D)&=1\\
        d_G(A,E)&=1\\
        d_G(B,C)&=2\\
                d_G(B,D)&=2\\
                d_G(B,E)&=1\\
                d_G(C,D)&=2\\
                d_G(C,E)&=1\\
                d_G(D,E)&=1.
        \end{split}
\]
and we get the same distance matrix $\Delta$ as for the above counterexample.
Here too the failure of conditional negativity stems from the lack of strict convexity. For this graph
\[
  d_G(A,B)+d_G(B,E)=  d_G(A,C)+d_G(C,E),
\]
with $C\not =E$. This will not happen for the graph on the right.
\section{A general theorem}
\setcounter{equation}{0}

\begin{theorem}
  \label{ththth}
  Let $\w$ be a continuous positive function on $\mathbb R$, with a zero set of measure $0$. Let $b\in(0,1]$.
  Let  $\mathcal V$ denote the vector space of continuous complex-valued functions on the real line such that $\int_{\mathbb R}|f(t)|^{b}\w(t)dt<\infty$.  Then,
  \begin{equation}
d(f,g)=\int_{\mathbb R}|f(t)-g(t)|^{b}\w(t)dt
\label{met}
\end{equation}
is a metric on $\mathcal V$, and
\[
D(f,g)=1-e^{-d(f,g)}
\]
is a normalized metric on $\mathcal V$, whose associated similarity metric is positive definite on $\mathcal V$.
\end{theorem}

\begin{proof} We proceed in a number of steps.\\

  STEP 1:{ \sl     For every real number $t$ the function
  \[
k_t(f,g)=|f(t)|^{b} +|g(t)|^{b}-|f(t)-g(t)|^{b}
\]
is positive definite on $\mathcal V$.}\smallskip

To prove Step 1, use Proposition \ref{hurst1} with $b$ instead of $2a$.\smallskip

STEP 2: {\sl The integral

  \[
\int_{\mathbb R}(|f(t)|^{b} +|g(t)|^{b}-|f(t)-g(t)|^{b})\w(t)dt,
\]
and its exponential
\[
  e^{\int_{\mathbb R}(|f(t)|^{b} +|g(t)|^{b}-|f(t)-g(t)|^{b})\w(t)dt}.
\]
are positive definite on $\mathcal V$.}\smallskip

The first claim follows from the fact that a sum of positive definite functions is positive definite, and approximating the integral. 
The second claim follows since the exponential of a positive definite function is positive definite.\\

STEP 3: {\sl The function
\[
  e^{-\int_{\mathbb R}(|f(t)-g(t)|^{b})\w(t)dt}
\]
is positive definite on $\mathcal V$.}\smallskip

This follows from
\[
    e^{-\int_{\mathbb R}|f(t)-g(t)|^{b}\w(t)dt}
=e^{-\int_{\mathbb R}|f(t)|^{b}\w(t)dt}\cdot
  e^{\int_{\mathbb R}(|f(t)|^{b} +|g(t)|^{b}-|f(t)-g(t)|^{b})\w(t)dt}\cdot e^{-\int_{\mathbb R}|g(t)|^{b}\w(t)dt}.
\]

Up to this stage we could have $b\in(0,2)$. To prove that $D(f,g)$ is a metric we prove that $d(f,g)$ is a metric, and we need $b\in(0,1]$ for that purpose.\smallskip

STEP 4: {\sl $d(f,g)$ is a metric on $\mathcal V$.}\smallskip

From Proposition \ref{hurst1} we see that $\delta(z,w)=|z-w|^b$ is a metric on $\mathbb C$ (see \eqref{xa}). Thus for $f,g,h\in\mathcal V$ and $t\in\mathbb R$,
\[
|f(t)-g(t)|^b\le |f(t)-h(t)|^b+|h(t)-g(t)|^b.
\]
Integrating over $\mathbb R$ with respect to the weight $\w$ we get that $d(f,g)$ satisfies the triangle inequality. The other properties for a metric are trivially verified. To conclude the proof of the theorem it remains to prove that $D(f,g)$ is indeed a metric, but this follows from the fact that $d$ is a metric,
by the already mentioned \cite[Lemma 12 p. 2370]{MR2522441}.
\end{proof}

\begin{remark}
  That $\delta(z,w)=|z-w|^b$ is a metric on $\mathbb C$ can also be seen from the fact that, for $b\in(0,1]$, $u^b$ is a metric when $u$ is. Furthermore,
  we note that \eqref{met} does not define a norm when $b\in(0,1)$.  
\end{remark}

We also note that one could replace $\w(t)dt$ by  a jump measure and obtain discrete counterparts of Theorem \ref{ththth}.
\section{More examples}
\setcounter{equation}{0}
For the first example, we note that if $d(x,y)$ is a distance on $X$ so is
\[
  \frac{d(x,y)}{1+d(x,y)}.
\]
The latter is a normalized metric.

\begin{example}
  We take $X=\mathbb C^n$ and $d(x,y)=\|x-y\|$. Then $D(x,y)=\frac{\|x-y\|}{1+|\|x-y\|}$ is a normalized metric, with associated similarity metric
  $S(x,y)=\frac{1}{1+\|x-y\|}$. The latter is positive definite on $\mathbb C^n$ as can be seen from the formula
  \[
\frac{1}{1+\|x-y\|}=\int_0^\infty e^{-u(1+\|x-y\|)}du.
    \]
\end{example}

\begin{example}
  We take $X=\mathbb C^n$. Then for $a\in(0,1)$ the function $d(x,y)=\|x-y\|^{a}$  is a metric. Restricting $x,y$ on the sphere of radius $(1/2)^{1/a}$, the associated similarity metric is
  \begin{equation}
    S(x,y)=1-\|x-y\|^{a}.    
  \end{equation}
  which is positive definite in the above sphere as the restriction of \eqref{hurst} to that sphere.
\end{example}

For the following example, see \cite{yianilos2002normalized}. A proof of the fact that one has a metric is also proved in \cite{alpay-it} using reproducing kernel methods.
  \begin{example} 
    We take $X=\mathbb C^n$ and $d(x,y)=\|x-y\|$. Then
\begin{equation}
  D(x,y)=\begin{cases}\, \dfrac{\|x-y\|}{\|x\|+\|y\|},\quad x,y\in\mathbb C^n, \, x\not=y,\\
      \,\, 0,\quad\hspace{14mm}\,\,\, x=y,
      \end{cases}
          \end{equation}
defines a metric, called the normalized Euclidean metric, with associated similarity metric
\[
  S(x,y)=\begin{cases}\,\,\dfrac{\|x\|+\|y\|-\|x-y\|}{\|x\|+\|y\|},\quad x\not=0\quad{\rm or}\quad y\not=0,\\
    \,\, 0,\,\hspace{2cm}\,\quad\hspace{1.4cm} x=y=0,\end{cases}
\]
which is positive definite on $\mathbb C^n$. The latter claim follows from the fact that the function $\|x\|+\|y\|-\|x-y\|$ (i.e. the function \eqref{hurst} with $a=1/2$)
is positive definite on $\mathbb C^n$, that
\[
\frac{1}{\|x\|+\|y\|}=\int_0^\infty e^{-t\|x\|}e^{-t\|y\|}dt
\]
is positive definite on $\mathbb C^n\setminus\left\{0\right\}$, and that the product of two positive definite functions is positive definite.
\end{example}
        \bibliographystyle{plain}
\def\cprime{$'$} \def\cprime{$'$} \def\cprime{$'$}
  \def\lfhook#1{\setbox0=\hbox{#1}{\ooalign{\hidewidth
  \lower1.5ex\hbox{'}\hidewidth\crcr\unhbox0}}} \def\cprime{$'$}
  \def\cprime{$'$} \def\cprime{$'$} \def\cprime{$'$} \def\cprime{$'$}
  \def\cprime{$'$}

\end{document}